\documentclass[oneside,english,reqno]{amsart}
\usepackage[T1]{fontenc}
\usepackage[latin9]{inputenc}
\usepackage{geometry}
\usepackage{color}
\usepackage{prettyref}
\usepackage{amsthm}
\usepackage{amsbsy}
\usepackage{amssymb}
\usepackage{graphicx}

\makeatletter

\newcommand{\noun}[1]{\textsc{#1}}

\numberwithin{equation}{section}
\numberwithin{figure}{section}
\numberwithin{table}{section}
\usepackage{enumitem}		
  \theoremstyle{remark}
  \newtheorem*{acknowledgement*}{\protect\acknowledgementname}
\theoremstyle{plain}
\newtheorem{thm}{\protect\theoremname}
  \theoremstyle{definition}
  \newtheorem{defn}[thm]{\protect\definitionname}

\usepackage{graphicx}
\usepackage{setspace}
\usepackage[hyphens]{url}
\usepackage[perpage,symbol]{footmisc}
\usepackage{multicol}
\usepackage{etoolbox}
\usepackage{bbm}
\usepackage{tensor}
\usepackage{enumitem}

\DefineFNsymbols*{lamportnostar}[math]{{(\dagger)}{(\ddagger)}{(\S)}{(\P)}{(\|)}{(\dagger\dagger)}{(\ddagger\ddagger)}}
\setfnsymbol{lamportnostar}

\newrefformat{prop}{Proposition~\ref{#1}}
\newrefformat{lem}{Lemma~\ref{#1}}
\newrefformat{sec}{§\ref{#1}}
\newrefformat{sub}{§\ref{#1}}
\theoremstyle{remark}
\newtheorem*{rems*}{Remarks}

\setitemize[0]{leftmargin=10pt,itemindent=0pt}
\setlist[enumerate]{leftmargin=*,widest=0}

\allowdisplaybreaks[2]

\renewenvironment{thebibliography}[1]{%
    \begin{oldthebibliography}{#1}%
      \setlength{\parskip}{0ex}%
      \setlength{\itemsep}{0ex}%
      \begin{small}
  }%
  {%
    \end{small}
    \end{oldthebibliography}%
  }

\makeatother

\usepackage{babel}
\usepackage{listings}
  \providecommand{\acknowledgementname}{Acknowledgement}
  \providecommand{\definitionname}{Definition}
\providecommand{\theoremname}{Theorem}

\begin{document}

\title[Hypermatrix Algebra]{A combinatorial approach to the algebra of hypermatrices }

\author{Edinah K. Gnang }
\begin{abstract}
We present two hypermatrix formulations of the Cayley\textendash Hamilton
theorem. One of the proposed formulation naturally extends to hypermatrices
the combinatorial interpretations of the classical Cayley\textendash Hamilton
theorem. We conclude by discussing an application of the theorem to
computing graph invariants which distinguish some non-isomorphic graphs
with isospectral adjacency matrices.
\end{abstract}

\date{\today}

\maketitle
\tableofcontents{}

\section{\label{sec:Introduction}Introduction.}

The importance of a graph theoretical perspective to the algebra of
matrices is well established\cite{Z,BC}. We show that insights provided
by a combinatorial lens on the algebra of matrices also shed light
on the algebra of multidimensional generalization of matrices called
hypermatrices. Formally, a hypermatrix denotes a finite set of numbers
whose distinct members are indexed by distinct elements of a Cartesian
product set of the form
\[
\left\{ 1,2,\cdots,n_{1}\right\} \times\left\{ 1,2,\cdots,n_{2}\right\} \times\cdots\times\left\{ 1,2,\cdots,n_{d}\right\} .
\]
Such a hypermatrix is said to be of order $d$ and of size $n_{1}\times n_{2}\times\cdots\times n_{d}$.
In particular matrices are second order hypermatrices. The algebra
of hypermatrices arises from attempts to extend to hypermatrices familiar
matrix algebra concepts \cite{MB94,GKZ,RK,GER}. A survey of important
hypermatix results can be found in \cite{Lim2013}. The discussion
here mostly focuses on the Bhattacharya-Mesner (BM) hypermatrix algebra
\cite{MB90,MB94}. On occasion we also discuss the general BM product
developed in \cite{GER,2014arXiv1411.6270G}. The general BM product
has the benefit of encompassing as special cases many other hypermatrix
products such as the Segre outer product, the contraction product
and the multilinear matrix multiplication described in detailed in
\cite{Lim2013}. Our main result are two new hypermatrix formulations
of the Cayley\textendash Hamilton theorem. The first of which extends
to hypermatrices combinatorial interpretations of the classical Cayley\textendash Hamilton
theorem described in \cite{BC,Z}, while the second formulation is
distinctively less combinatorial and more algebraic. The second formulation
has the benefit of bearing a close resemblance to the classical Cayley\textendash Hamilton
theorem. It also lends itself more easily to the computation of invariants.
Finally we discuss an application of the hypermatrix formulations
of the Cayley\textendash Hamilton theorem to computing graph invariants
which distinguish some non-isomorphic graphs whose adjacency matrices
are isospectral.
\begin{acknowledgement*}
We would like to thank Andrei Gabrielov for providing guidance and
inspiration while preparing this manuscript. We would like to thank
Vladimir Retakh and Ahmed Elgammal for patiently introducing us to
the theory of hypermatrices. We are grateful to Doron Zeilberger,
Ha Luu and Sowmya Srinivasan for helpful discussions and suggestions.
The author was supported by the the National Science Foundation, and
is grateful for the hospitality of the Institute for Advanced Study.
\end{acknowledgement*}

\section{\label{sec:Complexes-and-buildings}Overview of the Bhattacharya-Mesner
algebra}

We recall here for convenience of the reader the basic elements of
the Bhattacharya-Mesner (BM) algebra proposed in \cite{MB90,MB94}
as a generalization of the algebra of matrices.
\begin{defn}
The Bhattacharya-Mesner \cite{MB90,MB94} algebra generalizes the
classical matrix product 
\[
\mathbf{B}=\mathbf{A}^{(1)}\cdot\mathbf{A}^{(2)}
\]
where $\mathbf{A}^{(1)}$, $\mathbf{A}^{(2)}$, $\mathbf{B}$ are
matrices of sizes $n_{1}\times k$, $k\times n_{2}$, $n_{1}\times n_{2}$,
respectively, 
\[
b_{i_{1},i_{2}}=\sum_{1\le\textcolor{red}{j}\le k}a_{i_{1},\textcolor{red}{j}}^{(1)}\,a_{\textcolor{red}{j},i_{2}}^{(2)},
\]
to an $m$-operand hypermatrix product noted
\[
\mathbf{B}=\mbox{Prod}\left(\mathbf{A}^{(1)},\,\cdots,\mathbf{A}^{(m)}\right),
\]
where $\mathbf{B}$ is an $n_{1}\times\cdots\times n_{m}$ hypermatrix,
for $i=1,\cdots,\left(m-1\right)$, $\mathbf{A}^{(i)}$ is a hypermatrix
whose size is obtained by replacing $n_{i+1}$ by $k$ in the dimensions
of the hypermatrix $\mathbf{B}$, and $\mathbf{A}^{(m)}$ is a $k\times n_{2}\times\cdots\times n_{m}$
hypermatrix, 
\[
b_{i_{1},\cdots,i_{m}}=\sum_{1\le\textcolor{red}{j}\le k}a_{i_{1},\textcolor{red}{j},i_{3},\cdots,i_{m}}^{(1)}\cdots\,a_{i_{1},\cdots,i_{t},\textcolor{red}{j},i_{t+2},\cdots,i_{m}}^{(t)}\cdots\,a_{\textcolor{red}{j},i_{2},\cdots,i_{m}}^{(m)}.
\]
In the particular case of third order hypermatrices, $\mathbf{A}^{(1)}$,
$\mathbf{A}^{(2)}$, $\mathbf{A}^{(3)}$ and $\mathbf{B}$ are hypermatrices
of sizes $n_{1}\times k\times n_{3}$, $n_{1}\times n_{2}\times k$,
$k\times n_{2}\times n_{3}$ and $n_{1}\times n_{2}\times n_{3}$
respectively, 
\[
b_{i_{1},i_{2},i_{3}}=\sum_{1\le\textcolor{red}{j}\le k}a_{i_{1},\textcolor{red}{j},i_{2}}^{(1)}\,a_{i_{1},i_{2},\textcolor{red}{j}}^{(2)}\,a_{\textcolor{red}{j},i_{1},i_{2}}^{(3)}.
\]
The general BM product was introduced in \cite{GER} and noted 
\[
\mathbf{C}=\mbox{Prod}_{\mathbf{B}}\left(\mathbf{A}^{(1)},\,\cdots,\mathbf{A}^{(m)}\right).
\]
The hypermatrix $\mathbf{C}$ is an $n_{1}\times\cdots\times n_{m}$
hypermatrix, while the dimensions of the hypermatrix $\mathbf{A}^{(i)}$
for $i=1,\cdots,\:m-1$ is obtained by replacing $n_{i+1}$ by $k$
in the dimensions of $\mathbf{C}$ and $\mathbf{A}^{(m)}$ is a hypermatrix
of size $k\times n_{2}\times\cdots\times n_{m}$ similarly to the
BM product. Crucially, the general BM product differs from the BM
product in the fact that the general product involves an additional
input hypermatrix. The additional product input hypermatrix $\mathbf{B}$
is called the background hypermatrix and as such $\mathbf{B}$ must
be a cubic $m$-th order hypermatrix having all of its sides of length
$k$, 
\[
c_{i_{1},\cdots,i_{m}}=\sum_{1\le\textcolor{red}{j_{1}},\textcolor{red}{j_{2}},\textcolor{red}{\cdots},\textcolor{red}{j_{m}}\le k}a_{i_{1},\textcolor{red}{j_{2}},i_{3},\cdots,i_{m}}^{(1)}\cdots a_{i_{1},\cdots,i_{t},\textcolor{red}{j_{t+1}},i_{t+2},\cdots,i_{m}}^{(t)}\cdots a_{\textcolor{red}{j_{1}},i_{2},\cdots,i_{m}}^{(m)}\,b_{\textcolor{red}{j_{1}},\textcolor{red}{j_{2}},\textcolor{red}{\cdots},\textcolor{red}{j_{m}}}.
\]
Note that the original BM product is recovered by setting $\mathbf{B}$
to the Kronecker delta hypermatrix (i.e. the hypermatrix whose nonzero
entries all equal one and are located at the entries whose indices
all have the same value, in particular Kronecker delta matrices are
identity matrices).
\end{defn}

\section{Hypermatrix formulation of the Cayley\textendash Hamilton theorem.}

The classical Cayley\textendash Hamilton theorem, establishes a tight
upper bound for the dimension of the span of consecutive powers of
a generic $n\times n$ matrix. While it is clear that the dimension
of the span of consecutive Hadamard powers of a generic $n\times n$
matrix is $n^{2}$, it is surprising that the dimension of the span
of consecutive powers of a generic $n\times n$ matrix is at most
$n$. Similarly, hypermatrix formulations of the Cayley\textendash Hamilton
theorem establish tight upper bounds on the dimension of span of hypermatrix
\emph{powers}. Hypermatrix powers correspond to compositions of BM
products.

\subsection{First formulation of the Cayley\textendash Hamilton theorem.}

The first formulation of the Cayley\textendash Hamilton theorem is
based on the BM product introduced in \cite{MB90,MB94}. Recall that
the BM algebra is non associative. Consequently, the number of distinct
compositions of product a cubic hypermatrix $\mathbf{A}$ is determined
by the Fuss-Catalan numbers\cite{Lin}. In particular, a third order
hypermatrix $\mathbf{A}$ admits the following three distinct fifth
degree composition of product. 
\[
\mbox{Prod}\left(\mathbf{\mathbf{A}},\mathbf{A},\mbox{Prod}\left(\mathbf{A},\mathbf{A},\mathbf{A}\right)\right),
\]
\[
\mbox{Prod}\left(\mathbf{\mathbf{A}},\mbox{Prod}\left(\mathbf{A},\mathbf{A},\mathbf{A}\right),\mathbf{A}\right),
\]
\[
\mbox{Prod}\left(\mbox{Prod}\left(\mathbf{A},\mathbf{A},\mathbf{A}\right),\mathbf{\mathbf{A}},\mathbf{A}\right).
\]
Note that the BM product noted $\mbox{Prod}\left(\mathbf{A},\mathbf{A},\mathbf{A}\right)$
corresponds to a third degree power. Furthermore third order hypermatrices
admit by construction no even degree powers.
\begin{thm}
\label{thm:Cayley=002013Hamilton 1} The dimension of the span of
the vector space of third order cubic hypermatrix powers is maximal,
that is equal to the number of hypermatrix entries.
\end{thm}
For notational convenience, we restrict the discussion to third order
hypermatrices, however the argument presented here naturally extends
to hypermatrices of arbitrary order. 
\begin{proof}
We first observe each row-column slices of the powers of a generic
third order hypermatrix $\mathbf{A}$, can be expressed as some matrix
polynomial of the corresponding row-column slice of $\mathbf{A}$.
Consequently the upper bound on the dimension of the span of powers
of cubic hypermatrices of order $d$ and of side length $n$ is a
fixed polynomial in $n$ noted $p_{d}\left(n\right)$. Furthermore
the third order BM product is ternary, the number of distinct powers
of degree $2k+1$ is determined by the recurrence formula 
\begin{equation}
c_{3}=1,\quad c_{2k+1}=\sum_{0<i,j,i+j<2k+1}c_{i}\,c_{j}\,c_{2k+1-\left(i+j\right)}.\label{eq:recurrence formula}
\end{equation}
The recurrence \ref{eq:recurrence formula} is a special case of the
Fuss-Catalan numbers\cite{Lin} and in this particular case given
by 
\[
c_{2n+1}=\frac{{3n \choose n}}{2n+1}.
\]
as easily verified via the WZ method \cite{PWZ}. Furthermore, it
is clear that $p_{3}\left(n\right)$ is a polynomial of degree at
most $3$. Consequently by the polynomial argument it suffices to
exhibit explicit constructions of four hypermatrices $\mathbf{A}_{0}$,
$\mathbf{A}_{1}$, $\mathbf{A}_{2}$ and $\mathbf{A}_{3}$ of size
$n_{0}\times n_{0}\times n_{0}$, $n_{1}\times n_{1}\times n_{1}$,
$n_{2}\times n_{2}\times n_{2}$ and $n_{3}\times n_{3}\times n_{3}$
respectively such that 
\[
1\le n_{0}<n_{1}<n_{2}<n_{3}
\]
and most importantly, the span of the powers has maximal dimension.\\
Let $n_{0}=1$ and $\mathbf{A}_{0}$ be the third order hypermatrix
expressed 
\[
\mathbf{A}_{0}=\left[a_{1,1,1}=1\right].
\]
Let $n_{1}=2$ and $\mathbf{A}_{1}$ be determined by it's row column
$2\times2$ matrix slices given by 
\[
\mathbf{A}_{1}\left[:,:,0\right]=\left(\begin{array}{cc}
1 & 1\\
-1 & 1
\end{array}\right),\quad\mathbf{A}_{1}\left[:,:,1\right]=\left(\begin{array}{cc}
1 & 1\\
1 & 1
\end{array}\right).
\]
Let $n_{2}=3$ and $\mathbf{A}_{2}$ be determined by it's row column
$3\times3$ matrix slices given by
\[
\mathbf{A}_{2}\left[\,:,\,:,\,0\right]=\left(\begin{array}{rrr}
-1 & -1 & 45\\
0 & -8 & -1\\
3 & -79 & 1
\end{array}\right),\,\mathbf{A}_{2}\left[\,:,\,:,\,1\right]=\left(\begin{array}{rrr}
-3 & -1 & 2\\
-49 & 10 & -3\\
-6 & 2 & -1
\end{array}\right)
\]
\[
\mathbf{A}_{2}\left[\,:,\,:,\,2\right]=\left(\begin{array}{rrr}
-1 & 2 & -1\\
-1 & -1 & 0\\
-1 & 0 & -1
\end{array}\right)
\]
Finally, let $n_{3}=4$ and $\mathbf{A}_{3}$ be determined by it's
row column $4\times4$ matrix slices given by 
\[
\mathbf{A}_{3}\left[\,:,\,:,\,0\right]=\left(\begin{array}{rrrr}
2 & 0 & 2 & -1\\
-3 & 1 & 1 & 2\\
2 & -1 & 1 & 6\\
-1 & -3 & 0 & 20
\end{array}\right),\,\mathbf{A}_{3}\left[\,:,\,:,\,1\right]=\left(\begin{array}{rrrr}
0 & 0 & -1 & 3\\
0 & -1 & -20 & -1\\
2 & 1 & 2 & -1\\
3 & -1 & 1 & 0
\end{array}\right)
\]
\[
\mathbf{A}_{3}\left[\,:,\,:,\,2\right]=\left(\begin{array}{rrrr}
1 & 1 & 0 & -3\\
0 & 1 & 0 & 1\\
6 & -1 & -1 & 0\\
-2 & -2 & -5 & 2
\end{array}\right),\,\mathbf{A}_{3}\left[\,:,\,:,\,3\right]=\left(\begin{array}{rrrr}
-7 & -2 & -1 & 11\\
-1 & -1 & 3 & 78\\
-3 & 3 & 0 & -1\\
9 & 0 & 0 & 2
\end{array}\right).
\]
One easily verifies for $\mathbf{A}_{0}$, $\mathbf{A}_{1}$, $\mathbf{A}_{2}$
and $\mathbf{A}_{3}$ that the dimension of the vector space spanned
by the powers is respectively $1^{3}$, $2^{3}$, $3^{3}$ and $4^{3}$
respectively. This concludes the proof. 
\end{proof}
\smallskip{}

Having established the maximality of the span, Cramer's rule is used
to express the rational functions of the hypermatrix entries associated
with the linear dependence between of $n^{3}+1$ powers.

\subsection{Second formulation of the Cayley\textendash Hamilton theorem.}

Recall that the matrix powers can be computed via a recurrence formula
with initial conditions 
\[
\left\{ \mathbf{A}^{[0]}=\boldsymbol{\Delta},\:\mathbf{A}^{[1]}=\mathbf{A}\right\} 
\]
where 
\[
\left[\boldsymbol{\Delta}\right]_{i,j}=\begin{cases}
\begin{array}{cc}
1 & \mbox{ if }i=j\\
0 & \mbox{otherwise}
\end{array}\end{cases},
\]
and recurrence formula given by 
\[
\begin{cases}
\begin{array}{ccc}
\mathbf{A}^{[k+2]} & = & \mbox{Prod}_{\mathbf{A}^{[k]}}\left(\mathbf{A},\mathbf{A}\right)\\
\mathbf{A}^{[k+3]} & = & \mbox{Prod}_{\mathbf{A}^{[k+1]}}\left(\mathbf{A},\mathbf{A}\right)
\end{array}\end{cases}.
\]
Consequently, the classical Cayley\textendash Hamilton theorem establishes
the existence of sequence of rational functions 
\[
\left\{ \alpha_{k}\left(a_{1,1},\cdots,a_{n,n}\right)\right\} _{0\le k<n}\subset\mathbb{Q}\left(a_{1,1},\cdots,a_{n,n}\right)
\]
such that 
\[
\mathbf{0}_{n\times n}=\mathbf{A}^{[n]}+\sum_{0\le k<n}\mathbf{A}^{[k]}\,\alpha_{k}\left(a_{1,1},\cdots,a_{n,n}\right).
\]
The second hypermatrix formulation of the Cayley\textendash Hamilton
theorem is also defined by the recurrence 
\[
\left\{ \mathbf{A}^{[0]}=\boldsymbol{\Delta},\:\mathbf{A}^{[1]}=\mathbf{A}\right\} 
\]
where 
\[
\left[\boldsymbol{\Delta}\right]_{i,j,k}=\begin{cases}
\begin{array}{cc}
1 & \mbox{ if }i=j=k\\
0 & \mbox{otherwise}
\end{array}\end{cases},
\]
and recurrence formula given by 
\[
\begin{cases}
\begin{array}{ccc}
\mathbf{A}^{[k+2]} & = & \mbox{Prod}_{\mathbf{A}^{[k]}}\left(\mathbf{A},\mathbf{A},\mathbf{A}\right)\\
\mathbf{A}^{[k+3]} & = & \mbox{Prod}_{\mathbf{A}^{[k+1]}}\left(\mathbf{A},\mathbf{A},\mathbf{A}\right)
\end{array}\end{cases}.
\]

\begin{thm}
\label{thm:Cayley=002013Hamilton 2} The dimension of the span of
the vector space of third order cubic hypermatrix powers in the sequence
is maximal, that is equal to the number of hypermatrix entries.
\end{thm}
The proof of the theorem is similar to the previous proof in that
we observe each row-column slices of the powers of a generic third
order hypermatrix $\mathbf{A}$, can be expressed as some matrix polynomial
of the corresponding row-column slice of $\mathbf{A}$. Consequently
the upper bound on the dimension of the span of powers of cubic hypermatrices
of order $d$ and of side length $n$ is a fixed polynomial in $n$
noted $p_{d}\left(n\right)$.
\begin{proof}
The proof is similar to the proof given in the first formulation.
We describe hypermatrices $\mathbf{A}_{0}$, $\mathbf{A}_{1}$, $\mathbf{A}_{2}$
and $\mathbf{A}_{3}$ of size $n_{0}\times n_{0}\times n_{0}$, $n_{1}\times n_{1}\times n_{1}$,
$n_{2}\times n_{2}\times n_{2}$ and $n_{3}\times n_{3}\times n_{3}$
respectively such that 
\[
1\le n_{0}<n_{1}<n_{2}<n_{3}.
\]
The powers, of hypermatrices  the span of the powers has maximal dimension.\\
Let $n_{0}=1$ and $\mathbf{A}_{0}$ be the third order hypermatrix
expressed 
\[
\mathbf{A}_{0}=\left[a_{1,1,1}=1\right].
\]
Let $n_{1}=2$ and $\mathbf{A}_{1}$ be determined by its row column
$2\times2$ matrix slices given by 
\[
\mathbf{A}_{1}\left[:,:,1\right]=\left(\begin{array}{rr}
-7 & -7\\
1 & -1
\end{array}\right),\quad\mathbf{A}_{1}\left[:,:,2\right]=\left(\begin{array}{rr}
2 & 4\\
1 & -2
\end{array}\right).
\]
Let $n_{2}=3$ and $\mathbf{A}_{2}$ be determined by its row column
$3\times3$ matrix slices given by
\[
\mathbf{A}_{2}\left[\,:,\,:,\,1\right]=\left(\begin{array}{rrr}
-1 & 18 & 0\\
-3 & 0 & 5\\
2 & -1 & 2
\end{array}\right),\,\mathbf{A}_{2}\left[\,:,\,:,\,2\right]=\left(\begin{array}{rrr}
0 & -5 & -2\\
-3 & 1 & 1\\
1 & -2 & 0
\end{array}\right)
\]
\[
\mathbf{A}_{2}\left[\,:,\,:,\,3\right]=\left(\begin{array}{rrr}
-1 & 0 & 1\\
-1 & 2 & -14\\
6 & -3 & 1
\end{array}\right)
\]
Finally, let $n_{3}=4$ and $\mathbf{A}_{3}$ be determined by its
row column $4\times4$ matrix slices given by 
\[
\mathbf{A}_{3}\left[\,:,\,:,\,1\right]=\left(\begin{array}{rrrr}
18 & 0 & 0 & 1\\
0 & 3 & -1 & -1\\
0 & 52 & 4 & 5\\
-1 & -1 & -4 & 0
\end{array}\right),\,\mathbf{A}_{3}\left[\,:,\,:,\,2\right]=\left(\begin{array}{rrrr}
1 & 0 & -2 & 8\\
-2 & 1 & 1 & 1\\
4 & -2 & 6 & -2\\
-1 & -1 & 1 & -2
\end{array}\right)
\]
\[
\mathbf{A}_{3}\left[\,:,\,:,\,3\right]=\left(\begin{array}{rrrr}
10 & -1 & 0 & -1\\
1 & 0 & 1 & 3\\
1 & -1 & 0 & 0\\
0 & 1 & 13 & -1
\end{array}\right),\,\mathbf{A}_{3}\left[\,:,\,:,\,4\right]=\left(\begin{array}{rrrr}
4 & 12 & 2 & 0\\
-1 & -1 & -3 & 1\\
-1 & 1 & 0 & 0\\
155 & -1 & 0 & 0
\end{array}\right).
\]
One easily verifies for $\mathbf{A}_{0}$, $\mathbf{A}_{1}$, $\mathbf{A}_{2}$
and $\mathbf{A}_{3}$ that the dimension of the vector space spanned
by the powers is respectively $1^{3}$, $2^{3}$, $3^{3}$ and $4^{3}$
respectively. This concludes the proof. 
\end{proof}

\section{A combinatorial interpretation of the hypermatrix Cayley-Hamilton
theorem.}

Let $\mathbf{A}$ denote an $m\times n\times p$ third order hypermatrix.
We associate with $\mathbf{A}$ a directed tripartite $3$-uniform
hypergraph $H\left(\mathbf{A}\right)$. The hypergraph $H\left(\mathbf{A}\right)$
has $m$ vertices in the first partition, $n$ vertices in the second
partition and $p$ vertices in the third partition. The vertices in
the first, second and third partition are respectively colored red,
green and blue. The vertex coloring scheme is designed to establish
a one to one correspondence between entries of $\mathbf{A}$ and (
red, green, blue ) triplets of vertices in $H\left(\mathbf{A}\right)$.
More precisely, the directed hyperedge spanning the $i$-th red vertex
noted $\mbox{R}_{\textcolor{red}{i}}$, the $j$-th green vertex noted
$\mbox{G}_{\textcolor{green}{j}}$ and the $k$-th blue vertex noted
$\mbox{B}_{\textcolor{blue}{k}}$, is associated with the $a_{\textcolor{red}{i},\textcolor{green}{j},\textcolor{blue}{k}}$
hypermatrix entry. In short we say that $a_{\textcolor{red}{i},\textcolor{green}{j},\textcolor{blue}{k}}$
is the weight of the $\left(\mbox{R}_{\textcolor{red}{i}},\mbox{G}_{\textcolor{green}{j}},\mbox{B}_{\textcolor{blue}{k}}\right)$
hyperedge of $H\left(\mathbf{A}\right)$. The proposed directed tripartite
hypergraph $H\left(\mathbf{A}\right)$ described here is a natural
extension of the König directed bipartite graph associated with matrices
described in \cite{BC}.

\subsection{Composing Hypergraphs. }

By analogy to the matrix case, the König directed hypergraph yields
a combinatorial interpretation of the BM product. The hypergraph composition
is defined by the following vertex ( and induced edge ) identification
scheme. Consider tripartite hypergrahs $H\left(\mathbf{A}^{(1)}\right)$,
$H\left(\mathbf{A}^{(2)}\right)$, $H\left(\mathbf{A}^{(3)}\right)$
respectively associated with the $m\times t\times p$ hypermatrix
$\mathbf{A}^{(1)}$, the $m\times n\times t$ hypermatrix $\mathbf{A}^{(2)}$
and the $t\times n\times p$ hypermatrix $\mathbf{A}^{(3)}$. Incidentally,
the number of red vertices of $H\left(\mathbf{A}^{(1)}\right)$ equals
the number of red vertices of $H\left(\mathbf{A}^{(2)}\right)$. Similarly
the number of green vertices of $H\left(\mathbf{A}^{(2)}\right)$
also corresponds to the number of green vertices of $H\left(\mathbf{A}^{(3)}\right)$.
Finally the number of blue vertices of $H\left(\mathbf{A}^{(1)}\right)$
equals the number of blue vertices of $H\left(\mathbf{A}^{(3)}\right)$.
The size constraints, express the size requirement for the BM product
of $\mathbf{A}^{(1)}$, $\mathbf{A}^{(2)}$, and $\mathbf{A}^{(3)}$.
The result of the composition is a directed tripartite hypergraph
associated with an $m\times n\times p$ hypermatrix. As suggested
by the pairwise size constraints relating the hypergraph pair $\left(H\left(\mathbf{A}^{(1)}\right),H\left(\mathbf{A}^{(2)}\right)\right)$
the red vertices of $H\left(\mathbf{A}^{(1)}\right)$ are identified
according to their label with the red vertices of $H\left(\mathbf{A}^{(2)}\right)$.
Similarly, following the pairwise size constraints relating the hypergraph
pair $\left(H\left(\mathbf{A}^{(2)}\right),H\left(\mathbf{A}^{(3)}\right)\right)$
the green vertices of the hypergraph $H\left(\mathbf{A}^{(2)}\right)$
are identified according to their label to the green vertices of the
hypergraph $H\left(\mathbf{A}^{(3)}\right)$. Finally, following the
pairwise constraints relating the hypergraph pair $\left(H\left(\mathbf{A}^{(1)}\right),H\left(\mathbf{A}^{(3)}\right)\right)$
the blue vertices of $H\left(\mathbf{A}^{(1)}\right)$ are identified
according to their label with the blue vertices of $H\left(\mathbf{A}^{(3)}\right)$.
The final step of the identification consists in identifying vertices
of different colors according to their labels. Namely remaining green
vertices of $H\left(\mathbf{A}^{(1)}\right)$, the blue vertices of
$H\left(\mathbf{A}^{(2)}\right)$ as well as the red vertices of $H\left(\mathbf{A}^{(3)}\right)$
are identified according to their label values. Note that the last
identification step results into $t$ vertices whose color is neither
red, nor green nor blue. We assign the white color to such vertices.
Consequently, the weight associated with the $\left(\mbox{R}_{\textcolor{red}{r}},\mbox{G}_{\textcolor{green}{g}},\mbox{B}_{\textcolor{blue}{b}}\right)$
triplet of the hypermatrix resulting from the composition is given
by the summing over the white vertices as follows 
\[
\mbox{Weight of the }\left(\mbox{R}_{\textcolor{red}{r}},\mbox{G}_{\textcolor{green}{g}},\mbox{B}_{\textcolor{blue}{b}}\right)\mbox{ triplet in the composition}=\sum_{1\le w\le t}a_{\textcolor{red}{r},w,\textcolor{blue}{b}}^{(1)}\,a_{\textcolor{red}{r},\textcolor{green}{g},w}^{(2)}\,a_{w,\textcolor{green}{g},\textcolor{blue}{b}}^{(3)},
\]
the weighting of the resulting vertices correspond precisely to the
Bhattacharya-Mesner product. It therefore follows from the proposed
construction that 
\[
H\left(\mbox{Prod}\left(\mathbf{A},\,\mathbf{B},\,\mathbf{C}\right)\right)=\mbox{Composition}\left(H\left(\mathbf{A}\right),\,H\left(\mathbf{B}\right),\,H\left(\mathbf{C}\right)\right).
\]
It may be noted that each term of the form $a_{\textcolor{red}{r},w,\textcolor{blue}{b}}^{(1)}\,a_{\textcolor{red}{r},\textcolor{green}{g},w}^{(2)}\,a_{w,\textcolor{green}{g},\textcolor{blue}{b}}^{(3)}$
in the sum can be thought off as describing a tetrahedron construction
which connects the faces $\left(\textcolor{red}{r},w,\textcolor{blue}{b}\right)$,
$\left(\textcolor{red}{r},\textcolor{green}{g},w\right)$ and $\left(w,\textcolor{green}{g},\textcolor{blue}{b}\right)$.
It is therefore legitimate to deduce from the proposed identification
scheme that the edges (or sides) of the triangular faces are also
being appropriately identified. In particular, given an $n\times n\times n$
hypermatrix $\mathbf{A}$ with binary entries the sum 
\begin{equation}
\sum_{0\le r<g<b<n}\left[\mbox{Prod}\left(\mathbf{A},\,\mathbf{A},\,\mathbf{A}\right)\right]_{\textcolor{red}{r},\textcolor{green}{g},\textcolor{blue}{b}}
\end{equation}
counts the number tetrahedron construction possible using the hyperedge
from $H\left(\mathbf{A}\right)$. In particular for some particular
choice of ordered triplet $\left(\textcolor{red}{r},\textcolor{green}{g},\textcolor{blue}{b}\right)$
such that $0\le\textcolor{red}{r}<\textcolor{green}{g}<\textcolor{blue}{b}<n$
the entry $\left[\mbox{Prod}\left(\mathbf{A},\,\mathbf{A},\,\mathbf{A}\right)\right]_{\textcolor{red}{r},\textcolor{green}{g},\textcolor{blue}{b}}$
counts the number of tetrahedron in $H\left(\mathbf{A}\right)$ which
admit the ordered hyperedge $\left(\textcolor{red}{r},\textcolor{green}{g},\textcolor{blue}{b}\right)$
as one of the faces the tetrahedron. Furthermore the sum 
\[
\left[\mbox{Prod}\left(\mbox{Prod}\left(\mathbf{A},\,\mathbf{A},\,\mathbf{A}\right),\mathbf{A},\mathbf{A}\right)\right]_{\textcolor{red}{r},\textcolor{green}{g},\textcolor{blue}{b}}=\sum_{w_{1}}\left(\sum_{w_{0}}a_{\textcolor{red}{r},w_{0},\textcolor{blue}{b}}\,a_{\textcolor{red}{r},w_{1},w_{0}}\,a_{w_{0},w_{1},\textcolor{blue}{b}}\right)\,a_{\textcolor{red}{r},\textcolor{green}{g},w_{1}}\,a_{w_{1},\textcolor{green}{g},\textcolor{blue}{b}}
\]
counts the number of tetrahedral simplicial complex which can be constructed
by gluing two tetrahedrons at a face whose labels are of the form
$\left(\textcolor{red}{r},w_{1},\textcolor{blue}{b}\right)$ as depicted
in figure 4.1 where the face $\left(\textcolor{red}{r},w_{1},\textcolor{blue}{b}\right)$
is colored blue.
\begin{figure}
\includegraphics[scale=0.25]{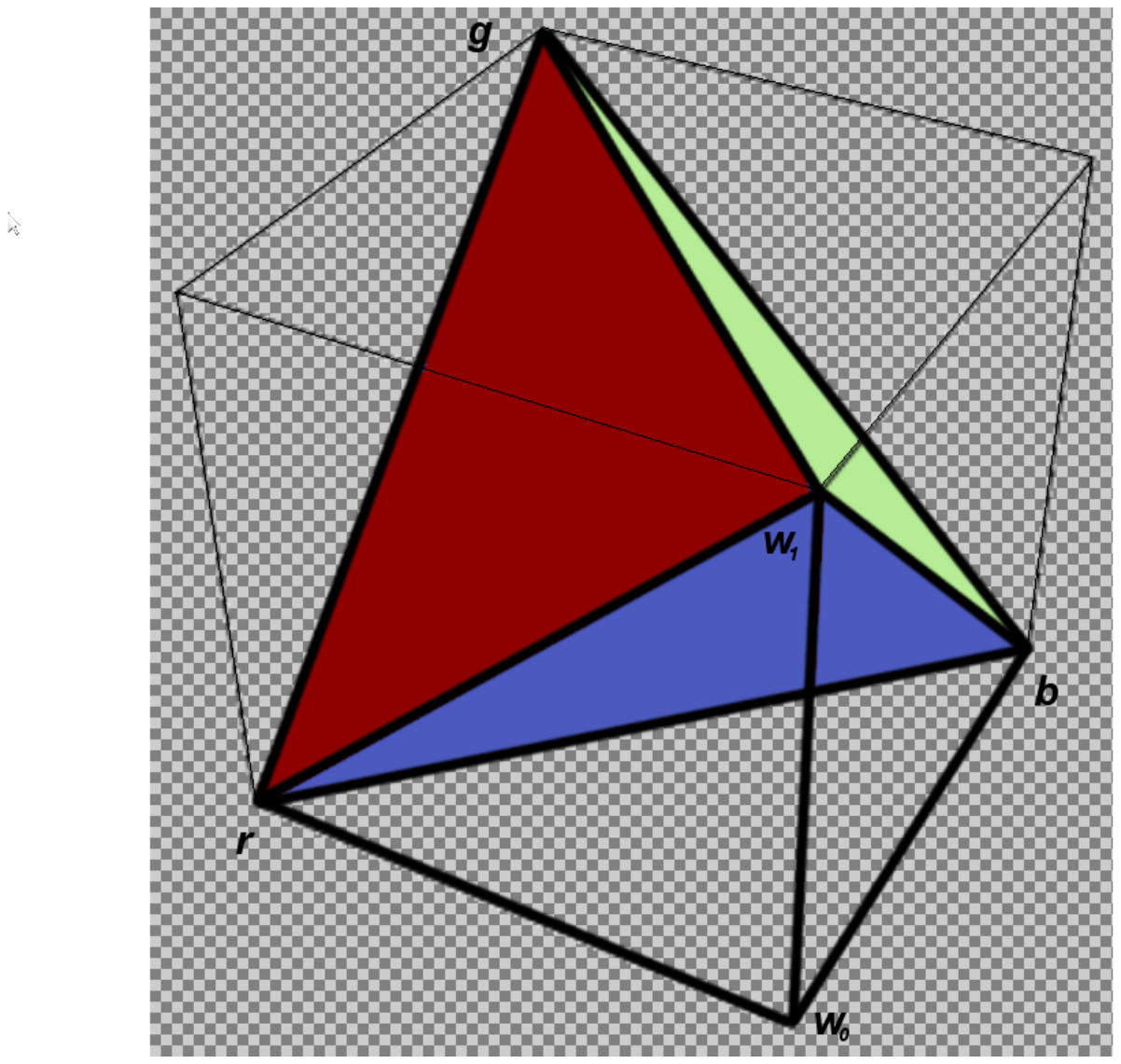}\centering\protect\caption{Gluing a tetrahedron on the face $\left(\textcolor{red}{r},w_{1},\textcolor{blue}{b}\right)$}

\end{figure}
 Furthermore the product 
\[
\left[\mbox{Prod}\left(\mathbf{A},\mbox{Prod}\left(\mathbf{A},\,\mathbf{A},\,\mathbf{A}\right),\mathbf{A}\right)\right]_{\textcolor{red}{r},\textcolor{green}{g},\textcolor{blue}{b}}=\sum_{w_{1}}a_{\textcolor{red}{r},w_{1},\textcolor{blue}{b}}\,\left(\sum_{w_{0}}a_{\textcolor{red}{r},w_{0},w_{1}}\,a_{\textcolor{red}{r},\textcolor{green}{g},w_{0}}\,a_{w_{0},\textcolor{green}{g},w_{1}}\right)\,a_{w_{1},\textcolor{green}{g},\textcolor{blue}{b}}
\]
counts the number of tetrahedral simplicial complex which can be constructed
by gluing two tetrahedrons at a face of whose labels are of the form
$\left(\textcolor{red}{r},\textcolor{green}{g},w_{1}\right)$ as depicted
in figure 4.2 where the face $\left(\textcolor{red}{r},\textcolor{green}{g},w_{1}\right)$
is colored red.
\begin{figure}
\includegraphics[scale=0.25]{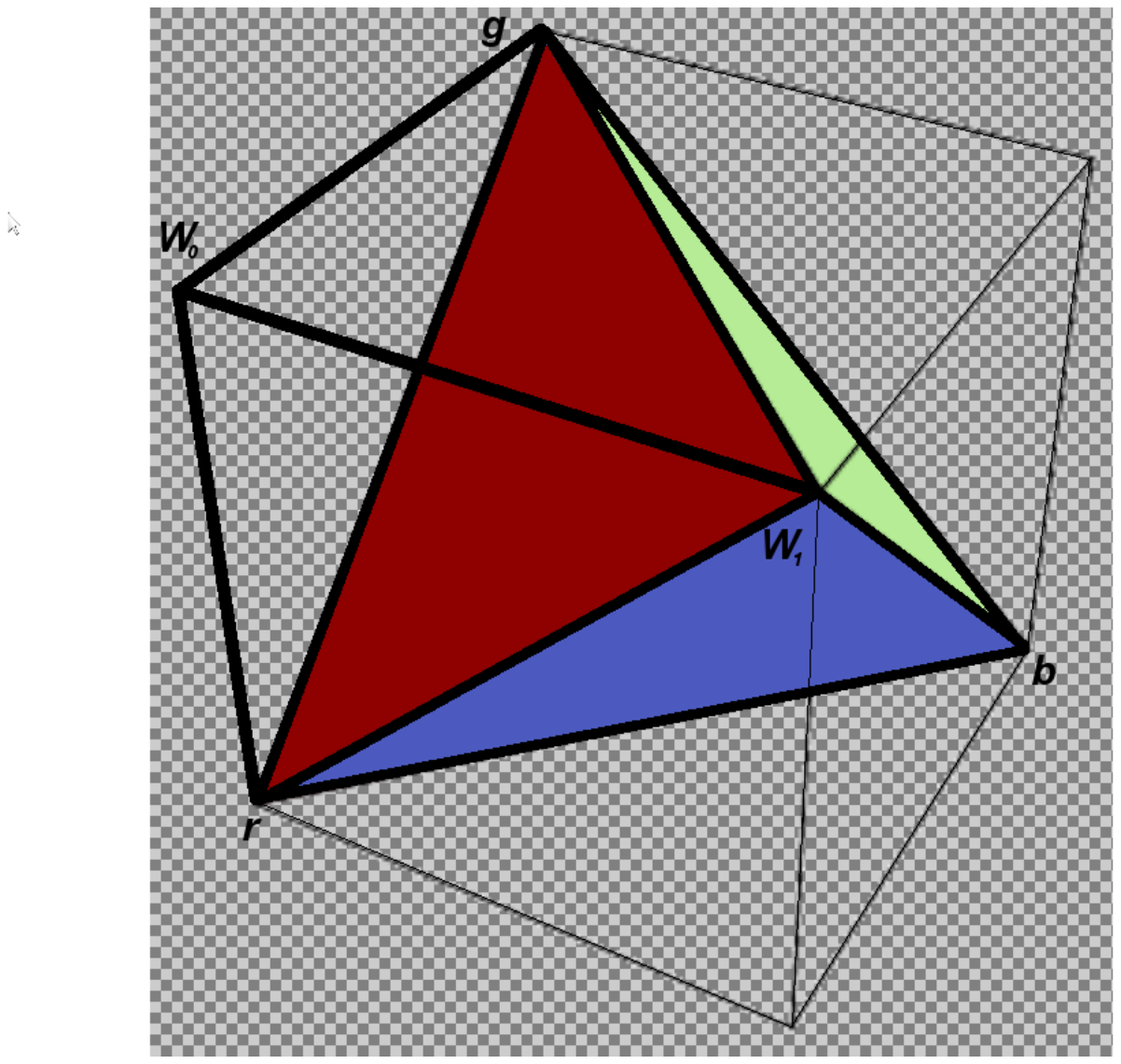}\centering\protect\caption{Gluing a tetrahedron on the face $\left(\textcolor{red}{r},\textcolor{green}{g},w_{1}\right)$ }
\end{figure}
 \\
Finally the product
\[
\left[\mbox{Prod}\left(\mathbf{A},\mathbf{A},\mbox{Prod}\left(\mathbf{A},\,\mathbf{A},\,\mathbf{A}\right)\right)\right]_{\textcolor{red}{r},\textcolor{green}{g},\textcolor{blue}{b}}=\sum_{w_{1}}a_{\textcolor{red}{r},w_{1},\textcolor{blue}{b}}\,a_{\textcolor{red}{r},\textcolor{green}{g},w_{1}}\,\left(\sum_{w_{0}}a_{w_{1},w_{0},\textcolor{blue}{b}}\,a_{w_{1},\textcolor{green}{g},w_{0}}\,a_{w_{0},\textcolor{green}{g},\textcolor{blue}{b}}\right)
\]
counts the number of tetrahedral simplicial complex which can be constructed
by gluing two tetrahedrons at a face whose labels are of the form
$\left(w_{1},\textcolor{green}{g},\textcolor{blue}{b}\right)$ as
depicted in figure 4.3 where the face $\left(w_{1},\textcolor{green}{g},\textcolor{blue}{b}\right)$
is colored green.
\begin{figure}
\includegraphics[scale=0.25]{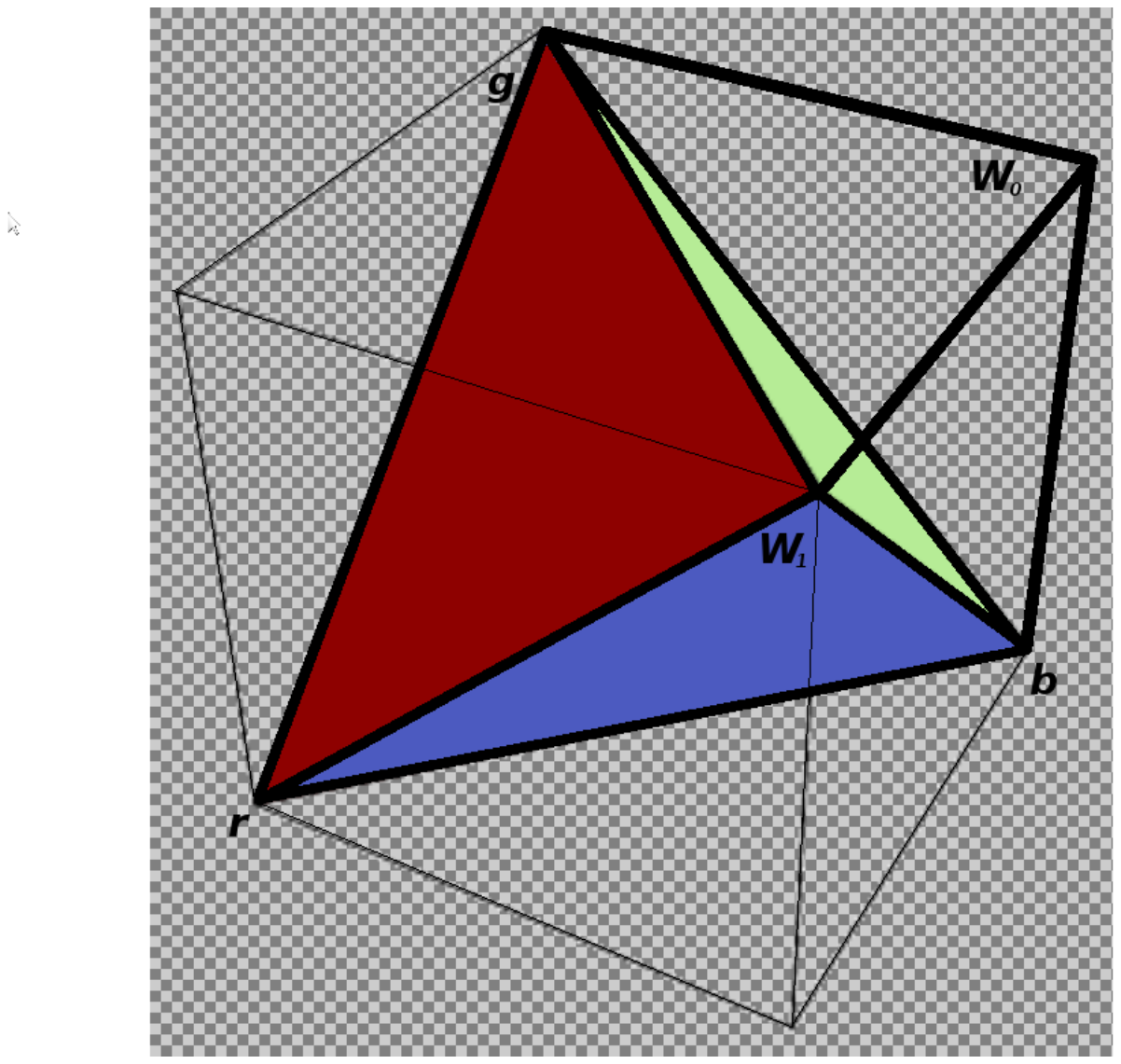}\centering\protect\caption{Gluing a tetrahedron on the face $\left(w_{1},\textcolor{green}{g},\textcolor{blue}{b}\right)$ }
\end{figure}

\section{Graph invariants via inflation.}

We shall aim to show here that the natural inflation scheme from graph
to hypergraphs introduced in \cite{AFW} combined with the combinatorial
invariants deduced from the generalization of the Cayley-Hamilton
theorem leads to symmetry breaking for some infinite families of cospectral
graphs. It is well known that the cospectrality for a pair of graphs
$G_{1}$ an $G_{2}$ is equivalent to the assertion that there exist
coefficients $\left\{ \alpha_{k}\right\} _{0<k\le n}$ such that
\[
\sum_{0<k\le n+1}\alpha_{k}\;\left(\#\mbox{ Walks of length }k\mbox{ connecting vertex }i\mbox{ to }j\mbox{ in }G_{1}\right)
\]
\[
=
\]
\[
\sum_{0<k\le n+1}\alpha_{k}\;\left(\#\mbox{ Walks of length }k\mbox{ connecting vertex }i\mbox{ to }j\mbox{ in }G_{2}\right)
\]
where $\alpha_{n+1}=1$, which algebraically expressed by the following
equality in terms of the adjacency matrices 
\begin{equation}
\left(\sum_{0<k\le n+1}\alpha_{k}\,\mathbf{B}^{k}\right)=0=\left(\sum_{0<k\le n+1}\alpha_{k}\,\mathbf{A}^{k}\right)
\end{equation}
Incidentally the property can be equivalently stated for an arbitrary
sequence of consecutive integer powers of $\mathbf{A}$, namely for
some arbitrary integer $\tau\ge0$ 
\begin{equation}
\left(\sum_{0<k\le n+1}\alpha_{k}\,\mathbf{B}^{\tau+k}\right)=0=\left(\sum_{0<k\le n+1}\alpha_{k}\,\mathbf{A}^{\tau+k}\right)
\end{equation}
This fact follows from the fact the vector space of powers of a matrix
has a span of dimension at most $n$ therefore we can more generally
state the cospectral invariance property by stating that 
\[
0=\sum_{0<k\le n+1}\left(\#\mbox{ Walks of length }\tau+k\mbox{ connecting }\left(i,j\right)\mbox{ in }G_{1}\right)\alpha_{k}
\]
\[
=
\]
\[
\sum_{0<k\le n+1}\left(\#\mbox{ Walks of length }\tau+k\mbox{ connecting }\left(i,j\right)\mbox{ in }G_{2}\right)\alpha_{k}
\]
 
\begin{thm}
\label{thm: Invariant}The sequence of sequence of Cayley\textendash Hamilton
coefficient are invariant under permutation of hypergaph vertices.
\end{thm}
The general argument of the proof is well illustrated for hypermatrices
of order $2$ and $4$ it will be immediately apparent how to extend
the argument to arbitrary even order hypermartices.
\begin{proof}
The proof that of invariance follows from the fact that the each BM
product corresponds to a sum over all vertices.
\end{proof}
Theorem\prettyref{thm: Invariant} establishes the Cayley\textendash Hamilton
coefficient as invariants hypermatrices. Similarly for hypermatrices
we may consider the equivalence classes between 3-uniform hypergraphs
induced by the 
\[
0=\sum_{0<k\le n^{3}+1}\alpha_{k}\:\left(\#\:k\mbox{\textendash Tetrahedral complex spanning }\left(\textcolor{red}{u},\textcolor{green}{v},\textcolor{blue}{w}\right)\mbox{ in }H_{1}\right)
\]
\[
=
\]
\[
\sum_{0<k\le n^{3}+1}\alpha_{k}\:\left(\#\:k\mbox{\textendash Tetrahedral complex spanning }\left(\textcolor{red}{u},\textcolor{green}{v},\textcolor{blue}{w}\right)\mbox{ in }H_{2}\right)
\]
(where a $k$-Tetrahrdral Simplex denotes a simplex using $k$ vertices
in addition to the boundary triangle vertices). The coefficient set
$\left\{ \alpha_{k}\right\} _{0<k\le n^{3}+1}$ where $\alpha_{n^{3}+1}=1$,
constitutes an invariant for hyperagraph under permutation the vertices
of the hypergraph. To show that such invariant are stronger then the
spectral invariant it suffices to consider the pair of adjacency matrices
with the smallest number of vertices which have the properties that
their adjacency matrices are cospectral. A tripartite 3-uniform hypergraph
is deduced from a graph as follows. We associate with to every directed
path of length two of the form $\textcolor{red}{v_{r}}\rightarrow\textcolor{green}{v_{g}}\rightarrow\textcolor{blue}{v_{b}}$,
an ordered hyperedge $\left(\mbox{R}_{\textcolor{red}{r}},\mbox{G}_{\textcolor{green}{g}},\mbox{B}_{\textcolor{blue}{b}}\right)$
of a hypergraph, thereby setting the $a_{rgb}$ entry of the adjacency
hypermatrix to $1$. We refer to such a construction as path adjacency
hypermatrix inflation. An easy rank argument on the compositions of
products reveals that the inflation scheme in conjunction with the
tetrahedral simplex counts indeed distinguishes the original two input
isospectral graphs and incidentally establishes the existence of an
infinite family of graphs for which the proposed inflation scheme
distinguishes isospectral non-isomorphic graphs.\bibliographystyle{amsalpha}
\bibliography{mybib}

\providecommand{\bysame}{\leavevmode\hbox to3em{\hrulefill}\thinspace}
\providecommand{\MR}{\relax\ifhmode\unskip\space\fi MR }
\providecommand{\MRhref}[2]{%
  \href{http://www.ams.org/mathscinet-getitem?mr=#1}{#2}
}
\providecommand{\href}[2]{#2}
\begin{thebibliography}{{Gna}14}

\bibitem[GER11]{GER}
E.~K. Gnang, A.~Elgammal, and V.~Retakh, \emph{A spectral theory for tensors},
  Annales de la faculte des sciences de Toulouse Mathematiques \textbf{20}
  (2011), no.~4, 801--841.

\bibitem[{Gna}14]{2014arXiv1411.6270G}
E.~K. {Gnang}, \emph{{Approximating the spectrum of matrices and
  hypermatrices}}, ArXiv e-prints (2014).

\bibitem[IGZ94]{GKZ}
M.M.~Kapranov I.M.~Gelfand and A.V. Zelevinsky, \emph{Discriminants, resultants
  and multidimensional determinant}, Birkhauser, Boston, 1994.

\bibitem[Ker08]{RK}
Richard Kerner, \emph{Ternary and non-associative structures}, International
  Journal of Geometric Methods in Modern Physics \textbf{5} (2008), 1265--1294.

\bibitem[Lim13]{Lim2013}
Lek-Heng Lim, \emph{Tensors and hypermatrices}, Handbook of Linear Algebra
  (Leslie Hogben, ed.), CRC Press, 2013.

\bibitem[Lin11]{Lin}
C-H Lin, \emph{Some combinatorial interpretations and applications of
  fuss-catalan numbers}, Discrete Mathematics \textbf{2011} (2011).

\bibitem[MB90]{MB90}
D.~M. Mesner and P.~Bhattacharya, \emph{Association schemes on triples and a
  ternary algebra}, Journal of combinatorial theory \textbf{A55} (1990),
  204--234.

\bibitem[MB94]{MB94}
D.~M. Mesner and P.~Bhattacharya, \emph{A ternary algebra arising from
  association schemes on triples}, Journal of algebra \textbf{164} (1994),
  595--613.

\bibitem[PWZ96]{PWZ}
Marko Petkov\v{s}ek, Herbert~S. Wilf, and Doron Zeilberger, \emph{A=b}, A.K.
  Peters, 1996.

\bibitem[RD08]{BC}
A.~Brualdi Richard and Cvetkovic Dragos, \emph{A combinatorial approach to
  matrix theory and its applications}, Chapman and Hall/CRC, 2008.

\bibitem[Zei85]{Z}
Doron Zeilberger, \emph{A combinatorial approach to matrix algebra}, Discrete
  Mathematics \textbf{56} (1985), 61--72.

\end{thebibliography}

\bigskip{}

\parbox[t]{1\columnwidth}{%
\noun{Department of Mathematics, Purdue University}\\
\noun{150 N. University Street, West Lafayette, IN 47907-2067}\\
E-mail:\texttt{ egnang@math.purdue.edu}%
}
\end{document}